\documentclass[10pt, a4paper]{amsart}

\usepackage{amsmath,amssymb,amscd,amsfonts}

\newtheorem{thm}{Theorem}[section]
\newtheorem{lem}[thm]{Lemma}
\newtheorem{rem}[thm]{Remark}
\newtheorem{prop}[thm]{Proposition}
\newtheorem{cor}[thm]{Corollary}

\newtheorem{claim}[thm]{Claim}

\newcommand\Pic{{\text{\rm Pic}}}
\newcommand\Ram{{\text{\rm Ram}}}
\newcommand\dimm{{\text{\rm dim}}}

\newcommand\PX{{\text {\rm Pic}}^0(X)}

\newcommand\DI{{\text {\rm dim}}   }

\newcommand\alb{\text{\rm a}}

\newcommand\Alb{\text{\rm A}}
\newcommand\lra{\longrightarrow}
\newcommand\ot{{\otimes}}
\newcommand\OO{{\mathcal{O}}}

\newcommand\PPP{{\mathbb{P}}}
\newcommand\ox{{\omega _X}}

\newcommand\FF{{\mathcal{F}}}

\newcommand\GG{{\mathcal{G}}}

\newcommand\Z{{\mathbb{Z}}}
\newcommand\Q{{\mathbb{Q}}}

\newcommand\III{{\mathcal{I}}}

\title[Varieties with $P_3(X)=4$ and $q(X)=\dim (X)$]
{Varieties with $P_3(X)=4$ and $q(X)=\dim (X)$ \footnote{2000 {
MSC. 14J10}}
}
\author {Alfred Jungkai Chen and Christopher D. Hacon }
\date{}
\begin{document}
\abstract We classify varieties with $P_3(X)=4$ and $q(X)=\dim (X)$.
\endabstract
\maketitle

\section{ Introduction}
Let $X$ be a smooth complex projective variety. When $\dim (X)\geq 3$
it is very hard to classify such varieties in terms of
their birational invariants.
Surprisingly, when $X$ has many holomorphic $1$-forms, it is sometimes
possible to achieve classification results in any dimension.
In \cite{Ka},
Kawamata showed that: {\it If $X$ is a smooth complex
projective variety with $\kappa (X)=0$ then the Albanese
morphism $\alb :X\lra \Alb (X)$ is surjective. If moreover,
$q(X)=\dim (X)$, then $X$ is
birational to an abelian variety.} Subsequently, Koll\'ar proved
an effective version of this result (cf. \cite{Ko2}):
{\em If $X$ is a smooth complex
projective variety with $P_m(X)=1$ for some $m\geq 4$, then the Albanese
morphism $\alb :X\lra \Alb (X)$ is surjective. If moreover,
$q(X)=\dim (X)$, then $X$ is
birational to an abelian variety.}
These results where further refined and expanded as follows:
\begin{thm} \label{T1}
{\bf (cf. \cite{CH1}, \cite{CH3}, \cite{HP}, \cite{Hac2})}
If $P_m(X)=1$ for some $m\geq 2$ or if $P_3(X)\leq 3$,
then the Albanese
morphism $\alb :X\lra \Alb (X)$ is surjective. If moreover $q(X)=\dim (X)$,
then:
\begin{enumerate}
\item If $P_m(X)=1$ for some $m\geq 2$, then $X$ is birational to
an abelian variety.
\item If $P_3(X)=2$, then $\kappa (X) =1$ and $X$ is a double cover
of its Albanese variety.
\item If $P_3(X)=3$, then $\kappa (X) =1$ and $X$ is a bi-double cover
of its Albanese variety.
\end{enumerate}
\end{thm}

In this paper we will prove a similar result for varieties with $P_3(X)=4$
and $q(X)=\dim (X)$. We start by considering the following examples:
\vskip.2cm
\noindent {\bf Example 1.}
Let $G$ be a group acting faithfully on a curve $C$ and
acting faithfully by translations on an abelian variety $\tilde{K}$, so that
$C/G=E$ is an elliptic curve and $\dim H^0(C, \omega _C^{\ot 3})^G=4$.
Let $G$ act diagonally on $\tilde{K}\times C$, then $X:= \tilde{K}\times C/G$
is a smooth projective variety with
$\kappa (X)=1$, $P_3(X)=4$ and $q(X)=\dim (X)$. We illustrate some
examples below:
\begin{enumerate}
\item $G=\Z _m$ with $m\geq 3$. Consider an elliptic curve $E$
with a line bundle $L$ of degree $1$.
Taking the normalization of the
$m$-th root of a divisor $B=(m-a)B_1+aB_2\in |mL|$
with $1\leq a\leq m-1$ and $m\geq 3$, one obtains a smooth curve $C$
and a morphism $g:C\lra E$ of degree $m$.
One has that $$g_* \omega _C=
\sum _{i=0}^{m-1} L^{(i)}$$
where $L^{(i)}=L^{\ot i}(-\lfloor
\frac{iB}{m}\rfloor )$ for $i=0,...,m-1$.

\item $G=\Z_2$. Let $L$ be a line bundle of degree $2$ over an elliptic curve
$E$. Let $C\lra E$ be the degree $2$ cover defined
by a reduced divisor $B\in |2L|$.

\item $G=(\Z _2 )^2$. Let $L_i$ for $i=1,2$ be line bundles of degree $1$
on an elliptic curve
$E$ and $C_i\lra E$ be degree $2$ covers
defined by disjoint reduced divisors $B_i\in |2L_i|$.
Then $C:=C_1\times _E C_2\lra E$ is a $G$ cover.

\item $G=(\Z _2 )^3$.  For $i=1,2,3,4$, let $P_i$ be distinct points
on an elliptic curve
$E$. For $j=1,2,3$ let
$L_j$ be line bundles of degree $1$ on $E$
such that $B_1=P_1+P_2\in |2L_1|$, $B_2=P_1+P_3\in |2L_2|$
and $B_3=P_1+P_4\in |2L_3|$.
Let $C_j\lra E$ be degree $2$ covers
defined by reduced divisors $B_j\in |2L_j|$.
Let $C$ be the normalization of
$C_1\times _E C_2\times _E C_3\lra E $, then $C$ is a $G$ cover.

\end{enumerate}
Note that (1) is ramified at 2 points.
Following \cite{Be} \S VI.12, one has that $P_2(X)=\dim H^0(
C, \omega _C ^{\ot 2})^G=2$ and $P_3(X)=\dim H^0(
C, \omega _C ^{\ot 3})^G=4$.
Similarly (2), (3), (4) are ramified along
4 points and hence $P_2(X)=P_3(X)=4$.
\vskip.2cm
\noindent {\bf Example 2.}
Let $q :A \lra S$ be a surjective morphism
with connected fibers from
an abelian variety of dimension $n\geq 3$ to an abelian surface.
Let $L$ be an ample line bundle on $S$ with $h^0(S,L)=1$,
$P\in \Pic ^0(A)$ with $P\notin \Pic ^0(S)$ and $P^{\ot 2}\in
\Pic ^0(S)$.
For $D$ an appropriate reduced divisor in $|L^{\ot 2}\ot
P^{\ot 2}|$, there is a degree $2$ cover $\alb :X\lra \Alb$ such that
$\alb _* (\OO _X)=\OO _\Alb \oplus (L\ot P)^\vee$.
One sees that $P_i(X)=1,4,4$ for $i=1,2,3$.

\vskip.2cm
\noindent {\bf Example 3.}
Let $q:\Alb \lra E_1\times E_2$ be a surjective morphism from an abelian
variety to the product of two elliptic curves, $p_i:\Alb \lra E_i$ the
corresponding morphisms, $L_i$ be line bundles of
degree $1$ on $E_i$ and $P,Q\in \Pic ^0(\Alb )$ such that
$P,Q$ generate a subgroup of $\Pic ^0(\Alb )/\Pic ^0(E_1\times E_2)$
which is isomorphic to $(\Z _2)^2$.
Then one has double covers $X_i\lra \Alb $ corresponding to divisors
$D_1\in |2(q_1^*L_1\ot P)|, D_2\in |2(q_2^*L_2\ot Q)|$.
The corresponding bi-double cover satisfies
$$\alb _* (\omega _X) = \OO _\Alb \oplus p_1^* L_1 \ot P
\oplus p_2^* L_2 \ot Q \oplus p_1^* L_1 \ot P  \ot p_2^* L_2 \ot Q$$
One sees that $P_i(X)=1,4,4$ for $i=1,2,3$.

\medskip

We will prove the following:
\begin{thm} \label{T2} Let $X$ be a smooth complex projective variety
with $P_3(X)=4$, then the Albanese morphism $\alb :X \lra \Alb$
is surjective (in particular $q(X)\leq \dim (X)$).
If moreover, $q(X)=\dim (X)$, then $\kappa (X)\leq 2$ and we have the
following cases:
\begin{enumerate}
\item
If $\kappa (X)=2$, then $X$ is birational either to a double cover
or to a bi-double cover of $\Alb$ as in Examples 2 and 3 and so
$P_2(X)=4$.
\item If $\kappa (X)=1$,
then $X$ is birational to the quotient $\tilde{K}\times C/G$ where
$C$ is a curve, $\tilde{K}$ is an abelian variety, $G$ acts faithfully on
$C$ and $\tilde{K}$. One has that either $P_2(X)=2$ and
$C\lra C/G$ is branched along 2 points with
inertia group $H \cong \Z _m$ with $m\geq 3$
or $P_2(X)=4$ and
$C\lra C/G$ is branched along 4 points with
inertia group $H \cong (\Z _2)^s$ with $s\in \{1,2,3\}$. See Example 1.
\end{enumerate}
\end{thm}

\bigskip

\noindent{\bf Acknowledgments.} The second author was partially supported
by NSA research grant no: MDA904-03-1-0101 and by a grant from the
Sloan Foundation.

\medskip
\noindent{\bf Notation and conventions.}
We work over the field of complex numbers.
We identify Cartier divisors and line bundles on
a smooth variety, and we use the additive and
multiplicative notation interchangeably. If
$X$ is a smooth projective variety, we let $K_X$ be a  canonical divisor,
so that
$\omega _X=\OO_X(K_X)$, and we denote by
$\kappa(X)$ the Kodaira dimension,  by
$q(X):=h^1(\OO_X)$ the {\em irregularity} and by $P_m(X):=h^0(\omega _X^{\ot m
})$
the {\em $m-$th
plurigenus}. We denote by  $\alb \colon X\to \Alb(X)$ the Albanese map  and by
$\Pic^0(X)$ the dual abelian variety to $\Alb (X)$ which parameterizes all
topologically trivial line bundles on $X$.
For a $\Q -$divisor $D$
we let $\lfloor D\rfloor$ be the integral part and
$\{D\}$ the fractional part. Numerical
equivalence is denoted by $\equiv$ and we write $D\prec E$ if
$E-D$ is an effective divisor. If $f\colon X\to Y$ is a
morphism, we write $K_{X/Y}:=K_X-f^*K_Y$ and we  often denote by  $F_{X/Y}$
the general fiber of
$f$. A $\Q$-Cartier divisor $L$ on a projective variety $X$ is nef
if for all curves $C\subset X$, one has $L.C\geq 0$.
For a surjective morphism of projective varieties
$f:X\lra Y$, we will say that a Cartier divisor $L$ on $X$ is $Y$-big
if for an ample line bundle $H$ on $Y$,
there exists a positive integer $m>0$ such that
$h^0(L^{\ot m}\ot f^*H ^\vee )>0$. The rest of the notation is
standard in algebraic geometry.

\section{Preliminaries }
\subsection{The Albanese map and the Iitaka fibration}
Let $X$ be a smooth projective variety. If $\kappa(X)>0$,
then the Iitaka fibration of $X$
is a morphism of projective varieties $f\colon X'\to Y$, with $X'$
birational to $X$ and
$Y$  of dimension
$\kappa(X)$, such that the general fiber of $f$ is smooth,
irreducible, of 
Kodaira dimension zero.  The Iitaka fibration is determined only up to
birational equivalence.
Since we are interested in questions of a birational nature,
we usually assume that
$X=X'$ and that $Y$ is smooth.

$X$ has {\em  maximal Albanese dimension} if
$\DI {(\alb _X(X))}=\DI (X)$.
We will need the following facts (cf.
\cite{HP} Propositions 2.1, 2.3, 2.12 and
Lemma 2.14 respectively).
\begin{prop}\label{albanese}Let $X$ be a smooth projective variety of maximal
Albanese dimension, and let
$f\colon X\to Y$ be the Iitaka fibration (assume $Y$ smooth). Denote by
$f_*\colon \Alb
(X)\to \Alb (Y)$ the homomorphism induced by
$f$ and  consider the commutative diagram:
$$
\CD
X &@>{\alb _X}>> & \Alb (X)\\
@V{f}VV& & @V{f_*}VV \\
Y &@>{\alb _Y}>> & \Alb (Y).\\
\endCD
$$

Then:
\begin{itemize}

\item[a)] $Y$ has maximal Albanese dimension;

\item[b)] $f_*$ is surjective and $\ker f_*$ is connected of  dimension
$\dim (X)-\kappa(X)$;

\item[c)] There exists an abelian variety
$P$ isogenous to $\ker f_*$ such that the general fiber of $f$ is birational to
$P$.
\end{itemize}
\end{prop}
Let $K:=\ker f_*$ and $F=F_{X/Y}$. Define
$$G:=ker \left( \Pic ^0(X)\to \Pic ^0(F)\right).$$
Then
\begin{lem}\label{LG}
$G$ is the union of finitely many translates of $\Pic ^0(Y)$
corresponding to the finite group $$\overline{G}:=G/\Pic ^0(Y)\cong ker \left(
\Pic ^0(K)\to \Pic ^0(F)\right).$$
\end{lem}

\subsection{Sheaves on abelian varieties}
Recall the following easy corollary of the theory of Fourier-Mukai transforms
cf. \cite{M}:
\begin{prop}\label{inclusion} Let $\psi\colon \FF \hookrightarrow
\GG$ be an inclusion of coherent sheaves on an abelian variety
$A$ inducing
isomorphisms $H^i(A,\FF \ot P)\to  H^i(A,\GG \ot P)$ for all
$i\geq 0$ and all $P\in \Pic ^0(A)$. Then $\psi$ is an isomorphism
of sheaves.
\end{prop}
Following \cite {M}, we will say that a coherent sheaf $\FF$ on an abelian
variety $A$ is I.T. 0 if $h^i(A,\FF \ot P)=0$ for all $i>0$.
We will say that an inclusion of coherent sheaves on $A$,
$\psi\colon \FF \hookrightarrow
\GG$ is an I.T. 0 isomorphism if $\FF,\GG$ are I.T. 0 and $h^0(\GG )=
h^0(\FF )$.  From the above proposition, it follows that every
I.T. 0 isomorphism $\FF \hookrightarrow
\GG$ is an isomorphism. We will need the following result:
\begin{lem}\label{L1}
Let $f:X\lra E$ be a morphism from a smooth projective variety to
an elliptic curve, such that $K_X$ is $E$-big. Then, for all
$P\in \Pic ^0(X)_{tors}$, $\eta \in \Pic ^0(E)$ and all 
$m\geq 2$, $f_* (\omega _X ^{\ot m}\ot P\ot f^* \eta )$ is I.T. 0.
In particular $$\deg (f_* (\omega _X ^{\ot m}\ot P\ot f^* \eta ))=
h^0(\omega _X ^{\ot m}\ot P\ot f^* \eta ).$$
\end{lem}
The proof of the above lemma is analogous to the proof of
Lemma 2.6 of \cite{Hac2}.
We just remark that it suffices to show that 
$f_* (\omega _X ^{\ot m}\ot P)$ is I.T. 0.
By \cite{Ko1}, one sees that $f_* (\omega _X ^{\ot m}\ot P)$
is torsion free and hence locally free on $E$.
By Riemann-Roch
$$h^0(\omega _X ^{\ot m}\ot P)=h^0(f_* (\omega _X ^{\ot m}\ot P))=\chi
(f_* (\omega _X ^{\ot m}\ot P))=\deg (f_* (\omega _X ^{\ot m}\ot P)).$$
\subsection{Cohomological support loci}
Let $\pi :X \lra A$ be a morphism from a smooth projective variety
to an abelian variety, $T\subset \Pic ^0(A)$ the translate of a subtorus
and $\FF$ a coherent sheaf on $X$.
One can define the cohomological support loci of $\FF$ as follows:
$$V^i(X,T,\FF ):= \{ P\in T|h^i(X,\FF \ot \pi ^* P)>0 \}.$$
If $T=\Pic ^0(X)$ we write
$V^i(\FF)$ or $V^i(X,\FF)$ instead of $V^i(X,\Pic ^0(X),\FF )$. When
$\FF = \omega _X$, the geometry of the loci
$V^i(\omega _X)$ is governed by the following result of Green and Lazarsfeld
(cf. \cite{GL}, \cite{EL}):
\begin{thm}\label{genvanish}{\rm \bf (Generic Vanishing Theorem)}
Let $X$ be a smooth projective variety.  Then:
\begin{itemize}
\item[a)]$V^i(\ox )$ has  codimension
$\ge i-(\dim (X)-\dim( \alb _X(X)))$;

\item[b)] Every irreducible component of $V^i(X,\ox)$ is
a translate of a sub-torus of\, $\PX$ by a torsion point
(the same also holds for the irreducible components of
$V^i_m(\ox ):=\{ P\in \Pic ^0(X)|h^i(X,\ox \ot P)\geq m\} $);

\item[c)] Let   $T$ be  an irreducible component of $V^i(\ox)$,  let
$P \in T$ be a point such that $V^i(\omega _X)$ is smooth at $P$, and let
$v
\in H^1(X, \OO
_X)\cong T_{P }\Pic ^0(X)$. If $v$ is not tangent to $T$, then the
sequence
$$H^{i-1}(X,\ox \otimes P) \stackrel {\cup v}  {\longrightarrow}
H^{i}(X,\ox \otimes P)  \stackrel {\cup v}  {\longrightarrow }
H^{i+1}(X, \ox \otimes P) $$
is exact. Moreover, if $P$ is a general point of $T$
and $v$ is
tangent to $T$ then
both maps   vanish;

\item[d)] If $X$ has maximal Albanese dimension, then there are inclusions:
$$ V^0(\ox )\supseteq V^1(\ox)
\supseteq \dots \supseteq V^n(\ox )=\{ \OO _X\}.$$

\item[e)] Let $f:Y\lra X$ be a surjective map of projective varieties, $Y$
smooth, then statements analogous to a), b), c) for $P\in \Pic ^0_{tors}
(Y)$ and d) above also hold for the sheaves
$R^if_*\omega _X$. More precisely we refer to \cite{CH3}, \cite {ClH} and
\cite{Hac5}.
\end{itemize}
\end{thm}
When $X$ is of maximal Albanese dimension,
its geometry is very closely connected to the properties of the loci
$V^i(\omega _X)$.
We recall the following two results from \cite{CH2}:
\begin{thm} \label{TCH2}
Let $X$ be a variety of maximal Albanese dimension. The translates
through the origin of the irreducible components of $V^0(\omega
_X)$ generate a subvariety of $\Pic ^0(X)$ of dimension $\kappa
(X) -\dim (X)+q(X)$. In particular, if $X$ is of general type then
$V^0(X,\omega _X )$ generates $\Pic ^0(X)$.
\end{thm}
\begin{prop}\label{PCH2}
Let $X$ be a variety of maximal Albanese dimension and $G,Y$
defined as in Proposition \ref{albanese}. Then

\begin{itemize}
\item[a)]$V^0(X,\Pic ^0 (X), \omega _X )\subset G$;

\item[b)] For every $P\in G$,
the loci $V^0(X, \Pic ^0 (X), \omega _X)\cap \left( P+\Pic ^0(Y)\right)$
are non-empty;

\item[c)] If $P$ is an isolated point of $V^0(X, \Pic ^0 (X), \omega _X)$,
then $P=\OO _X$.
\end{itemize}
\end{prop}
The following result governs the geometry of $V^0(\omega _X ^{\ot m})$
for all $m\geq 2$:
\begin{prop}\label{Pm} Let $X$ be a smooth projective variety
of maximal Albanese dimension,
$f\colon X\to Y$ the Iitaka fibration (assume $Y$ smooth) and
$G$ defined as in Proposition \ref{albanese}.
If $m\ge 2$, then $V^0(\omega _X ^{\ot m})=G$.
Moreover, for any fixed $Q\in V^0(\omega _X ^{\ot m})$, and all
$P\in \Pic ^0(Y)$ one has $h^0(\omega _X^{\ot m}\ot   Q
\ot P)=h^0(\omega_X^{\ot m}\ot Q).$
\end{prop}
We will also need the following lemma proved in \cite{CH2} \S 3.
\begin{lem} \label{L7} Let $X$ be a smooth projective variety
and $E$ an effective $\alb _X$-exceptional divisor on $X$.
If $\OO _X (E)\ot P$ is effective for some $P\in \Pic ^0(X)$,
then $P=\OO _X$.
\end{lem}
The following result is due to
Ein and Lazarsfeld (see \cite{HP} Lemma 2.13):
\begin{lem}\label{Lel} Let $X$ be a variety such 
that $\chi (\omega _X )=0$
and such that $\alb _X:X \lra \Alb (X)$ is surjective and generically finite.
Let $T$ be an irreducible component of $V^0(\omega _X)$,
and let $\pi _E :X \lra E:=\Pic ^0(T)$ be the morphism induced by the
map $\Alb (X)\lra \Pic^0
(\Pic ^0(X))\lra E$ corresponding to the inclusion
$T\hookrightarrow \Pic ^0(X)$.

Then there exists a divisor $D_T\prec R:=\Ram (\alb _X )=K_X$,
vertical with respect to $\pi _E$ (i.e. $\pi _E(D_T)\ne E$),
such that for general
$P\in T$, $G_T:=R-D_T$ is a fixed
divisor of each of the linear series $|K_X+P|$.
\end{lem}
We have the following useful Corollary
\begin{cor}
\label{C9} In the notation of Lemma \ref{Lel}, if $\dim (T)=1$, then
for any $P\in T$, there exists a line bundle of degree $1$ on
$E$ such that $\pi _E ^* L_P\prec K_X +P$.
\end{cor}
\begin{proof} By \cite{HP} Step 8 of the proof of Theorem 6.1, for general
$Q\in T$, there exists a line bundle of degree 1 on
$E$ such that $\pi _E ^* L_Q\prec K_X +Q$. Write $P=Q+\pi ^*\eta$ where
$\eta \in \Pic ^0(E)$. Then, since 
$$h^0(\omega _X \ot P\ot \pi ^* (L_Q\ot \eta )^\vee) = 
h^0(\pi _*(\omega _X\ot Q)\ot L_Q^\vee )\ne 0,$$ one sees that there is
an inclusion $\pi ^* (L_Q\ot \eta )\lra \omega _X \ot P$.
\end{proof}
Recall the following result (cf. \cite{Hac2} Lemma 2.17):
\begin{lem}\label{claimA}Let $X$ be a smooth projective variety,
let $L$ and $M$ be line bundles on $X$, and let  $T\subset \Pic
^0(X)$ be an irreducible subvariety of dimension $t$. If for all
$P\in T$, $\dim |L+ P|\geq a$ and $\dim |M-P|\geq b$, then $\dim
|L+ M|\geq a+b+t$.
\end{lem}
\begin{lem} \label{fiber} Let $T$ be a $1$-dimensional component of
$V^0(\omega _X)$, $E:=T^\vee$ and $\pi :X\lra E$ the induced
morphism. 
Then $P|F\cong \OO _F$ for all $P\in T$.
\end{lem}
\begin{proof} Let $G_T,D_T$ be as in Lemma \ref{Lel},
then for $P\in T$ we have
$|K_X+P|=G_T+|D_T+P|$ and hence the divisor $D_T+P$ is effective.
It follows that $(D_T+P)|F$ is also effective.
However $D_T$ is vertical with respect to $\pi$ and
hence $D_T|F\cong \OO _F$. By Lemma \ref{L7}, one sees that
$P|F \cong \OO _F$.
\end{proof}

\section{Kodaira dimension of Varieties with $P_3(X)=4,q(X)=\dimm (X)$ }
The purpose of this section is to study the Albanese map and
Iitaka fibration of varieties with $P_3=4$ and $q=\dimm (X)$. 
We will show that: 1) the Albanese map is surjective,
2) the image of the Iitaka fibration is an abelian
variety (and hence the Iitaka fibration factors through the Albanese map),
3) we have
that $\kappa(X) \le 2$.

We begin by fixing some notation. We write $$V_0(X,
\omega_X)=\cup_{i \in I} S_i$$ where $S_i$ are irreducible
components. Let $T_i$ denote the translate of $S_i$ 
passing through the origin and $\delta_i:=\dimm (S_i)$.  In
particular,  $S_0$ denotes the component contains the origin. For
any $i,j \in I$, let $\delta_{i,j}:=\dimm (T_i \cap T_j)$.

Recall that $V_0(X,\omega_X) \subset G \to \bar{G}:=G/\Pic^0(Y)$.
For any $\eta \in \bar{G}$, let $S_\eta$ denote a maximal
dimensional component which maps to $\eta$.
If $X$ is of maximal Albanese dimension with $q(X)=
\dimm (X)$, then its Iitaka fibration image $Y$ is of maixmal
Albanese dimension with $q(Y)=\dimm (Y) = \kappa(X)$. Moreover, by
Proposition \ref{PCH2}, one has $ \delta_i \ge 1, \forall i \ne
0$.

Now let $Q_i$ ($Q_\eta$ resp.) be a general torsion element in
$S_i$ ($S_\eta$ resp.), we denote by $P_{m,i}:=h^0(X, \omega_X
^{\ot m}\ot Q_i)$
($P_{m,\eta}$ resp.). Proposition
\ref{Pm} can be rephrased as

\begin{equation}\label{pm}
P_{m,\eta}=P_{m,\eta +\zeta}\ \ \ \forall \eta \in \Pic ^0(X),\ 
\zeta \in \Pic ^0(Y),\ m\geq 2 .
\end{equation}

By Lemma \ref{claimA} 
one has, for any $\eta, \zeta \in \bar{G}$,

\begin{equation}\label{pluri}
\left\{
\begin{array}{l}
P_{2,\eta+\zeta} \ge
P_{1,\eta}+P_{1,\zeta} +
\delta_{\eta,\zeta}-1,\\
 P_{2, 2 \eta} \ge 2 P_{1,\eta} + \delta_\eta -1,\\
P_{3,\eta+\zeta} \ge P_{1,\eta}+P_{2,\zeta} +\delta_\eta -1.
\end{array} \right.
\end{equation}

The following lemma is very useful when $\kappa \ge 2$.
\begin{lem}\label{elliptic}
Let $X$ be a variety of maximal Albanese dimension
with $\kappa (X) \ge 2$. Suppose that there is a
surjective morphism $\pi: X \to E$ to an elliptic curve $E$, and
suppose that there is an inclusion $\varphi:\pi^* L \to \omega _X^{\ot m}
\ot P$ for some $m \ge 2$, $P_{|F}=\OO_F$ where $F$ is a general fiber
of $\pi$ and $L$ is an ample line bundle
on $E$. Then the
induced map $L \to \pi_* (\omega_X^{\ot m} \ot P)$ is not an isomorphism,
$rank (\pi_* (\omega_X^{\ot m} 
\ot P)) \ge 2$ and $h^0(X,\omega_X^{\ot m} \ot P) >
h^0(E,L)$.

\end{lem}
\begin{proof}
By the easy addition theorem, $\kappa(F) \ge 1$. Hence by Theorem
1.1, $P_m(F) \ge 2$ for $m \ge 2$. The sheaf 
$\pi_* (\omega_X^{\ot m} \ot P)$ has rank equal
to $h^0(F,\omega_X^{\ot m} \ot P |_F)=h^0(F, \omega_F^{\ot m}) \ge 2$.
Therefore, $L \to \pi_* (\omega_X^{\ot m} \ot P)$ is not an isomorphism.
Since they are non-isomorphic I.T.0 sheaves, 
it follows that $h^0(\pi_* (\omega_X^{\ot m} \ot
P)) > h^0(L)$.
\end{proof}

\begin{cor}\label{ell}
Keep the notation as in Lemma \ref{elliptic}. If there is a
morphism $\pi': X \to E'$ and an inclusion $\pi'^* L'
\hookrightarrow \omega_X \ot P^\vee$ for some ample $L'$ on $E'$ and
$P\in \Pic ^0(X)$ with $P_{|F'}=\OO_{F'}$, then for all $m\geq 2$
$$ P_{m+1}(X) \ge 2+ h^0(X,\omega_X^{\ot m} \ot P) > 2+h^0(E',L').$$
\end{cor}

\begin{proof}
The inclusion $\pi'^* L' \hookrightarrow \omega_X \ot P^\vee$
induces an inclusion $$\pi'^* L' \ot \omega_X^{\ot m} \ot P
\hookrightarrow \omega_X^{\ot m+1}.$$ 
By Riemann-Roch, one has
$$ P_{m+1}(X) \ge h^0(E', L' \ot \pi'_* (\omega_X^{\ot m} \ot
P))\geq h^0(E',\pi'_* (\omega_X^{\ot m} \ot P))+rank(\pi'_* (\omega_X^
{\ot m} \ot
P)).$$ Proposition \ref{PCH2}, there exists $\eta \in \Pic ^0(Y)$ such that
$h^0(\omega _X^{\ot m-1}\ot P^{\ot 2}\ot \eta)\ne 0$ and
hence there is an inclusion
$$\pi'^* L' \hookrightarrow \omega _X ^{\ot m}\ot P\ot \eta.$$
By Proposition \ref{Pm} and Lemma \ref{elliptic},
$$h^0(X, \omega_X^{\ot m} \ot P)= h^0(X, \omega_X^{\ot m} \ot P\ot \eta )>
h^0(E',L').$$
\end{proof}

\begin{rem}\label{1dim}
Let $X$ be a variety with $\kappa(X) \ge 2$. Suppose that there
is a 1-dimensional component $S_i \subset V^0(\omega_X)$. We often
consider the induced map $\pi: X \to E:=T_i^\vee$. It is 
easy to see that $\pi$ factors through the Iitaka fibration. By Corollary
\ref{C9} and Lemma \ref{fiber}, there is an inclusion
$\varphi:\pi^* L \to \omega _X\ot P$ for some $P\in \Pic ^0(X)$ with
$P_{|F}=\OO_F$ and some ample
line bundle $L$ on $E$. In what follows, we will often apply
Lemma \ref{elliptic} and Corollary \ref{ell}
to this situation.
\end{rem}
\begin{lem}\label{P2} Let $X$ be a variety of maximal Albanese dimension
with $\kappa (X) \ge 2$ and $P_3(X)=4$. 
Then for any $\zeta \in G-\Pic ^0(Y)$,
one has $P_{2,\zeta}\leq 2$.
\end{lem}
\begin{proof} If $P_{2,\zeta}\geq 3$, then by 
(\ref{pluri}) and Proposition \ref{PCH2}, one sees that 
$V^0(\omega _X )\cap (\Pic ^0(Y)-\zeta )$ consists of $1$-dimensional
components. Let $S$ be one such component and $\pi :X\lra E:=S^\vee$
be the induced morphism. Then there is an ample line bundle $L$ on 
the elliptic curve $E$ and an inclusion $L\lra \pi _*(\omega _X \ot Q)$ 
for some $Q\in \Pic ^0(Y)-\zeta$. By Corollary \ref{ell},
$P_3(X)\geq 2+P_{2,\zeta}\geq 5$ which is impossible.
\end{proof}

\begin{thm}\label{surj}
Let $X$ be a smooth projective variety with $P_3(X)=4$, then the
Albanese morphism $\alb :X \lra \Alb$ is surjective.
\end{thm}
\begin{proof} We follow the proof of Theorem 5.1 of \cite{HP}.
Assume that $\alb :X \lra \Alb$ is not surjective, then we may assume
that there is a morphism $f:X\lra Z$ where $Z$ is a smooth variety of
general type, of dimension at least $1$, such that its Albanese map
$\alb _Z:Z\lra S$ is birational onto its image. By the proof of
Theorem 5.1 of \cite{HP}, it suffices to consider the cases in which
$P_1(Z)\leq 3$ and hence $\dim (Z)\leq 2$.
If $\dim (Z)=2$, then $q(Z)=\dim (S)\geq 3$ and since $\chi (\omega _Z)>0$,
one sees that $V^0(\omega _Z)=\Pic ^0(S)$.
By the proof of Theorem 5.1 of \cite{HP}, one has that for generic
$P\in \Pic ^0(S)$,
$$P_3(X)\geq h^0(\omega _Z\ot P) +h^0(\omega _X^{\ot 3}\ot
f^*\omega _Z^\vee \ot P)+\dim (S) -1\geq 1+2+3-1\geq 5.$$
This is a contradiction, so we may assume that $\dim (Z)=1$.
It follows that $g(Z)=q(Z)=P_1(Z)\geq 2$ and one may write
$\omega _Z=L^{\ot 2}$ for some ample line bundle $L$ on $Z$.
Therefore, for general $P\in \Pic ^0(Z)$, one has that $h^0(\omega
_Z\ot L\ot P)\geq
2$ and proceeding as in the proof of
Theorem 5.1 of \cite{HP}, that $h^0(\omega _X^{\ot 3}\ot f^*(\omega _Z
\ot L)^\vee \ot P)\geq 2$.
It follows as above that
$$P_3(X)\geq h^0(\omega _Z\ot L\ot P) +h^0(\omega _X^{\ot 3}\ot
f^*(\omega _Z\ot L)^\vee \ot P)+\dim (S) -1\geq
2+2+2-1\geq 5.$$
This is a contradiction and so $\alb :X \lra \Alb$ is surjective.
\end{proof}

\begin{prop} \label{gt} Let $X$ be a smooth projective variety 
with $P_3(X)=4,q(X)=\dim (X)$, then  
\begin{enumerate}
\item $X$ is not of general type and
\item if $\kappa (X)\geq 2$, then
$$V^0(\omega _X )\cap f^* \Pic ^0(Y)=\{\OO _X \}.$$
\end{enumerate}
\end{prop}
\begin{proof}
If $\kappa(X)=1$, then clearly $X$ is not of general type as
otherwise $X$ is a curve with
$P_3(X)=5g-5 > 4$. We thus assume that $\kappa(X) \ge
2$. It suffices to prove (2) as then (1) will follow from
Theorem \ref{TCH2}.

If all points of $V^0(\omega _X )\cap f^* \Pic ^0(Y)$ are
isolated, then the above statement follows from Proposition
\ref{PCH2}. Therefore, it suffices to prove that $\delta_0=0$.
(Recall that $\delta_0$ is the maximal dimension of a component in
$\Pic^0(Y)$.)

Suppose that $ \delta_0 \ge 2$. Then by (\ref{pluri}) and
Proposition \ref{Pm}, one has
$$ P_2 \ge 1+1+\delta_0 -1 \ge 3,\ \ \ \ \ \ \ 
 P_3 \ge 3+1+\delta_0 -1 \ge 5$$
which is impossible.

Suppose now that $ \delta_0 =1$, i.e. there is a $1$-dimensional
component $T\subset V^0(\omega _X)\cap f^*\Pic ^0(Y)$. Let $\pi
:X\lra E:=T^\vee$ be the induced morphism. By Corollary \ref{C9},
for some general $P\in T$, there exists a line bundle of degree 1
on $E$ and an inclusion $\pi ^*L\lra \omega _X \ot P$. By Lemma
\ref{fiber}, $P|F_{X/E}\cong \OO _{F_{X/E}}$.

We consider the inclusion $\varphi: L^{\ot 2}\lra \pi
_*(\omega_X^{\ot 2}\ot P^{\ot 2})$. By Lemma \ref{elliptic},
 one
sees that $h^0(\omega_X^{\ot 2}\ot P^{\ot 2})\geq 3$, and
$rank (\pi_*(\omega_X^{\ot 2} \ot P^{\ot 2})) \ge 2$. So
$$P_3(X)=h^0(\omega _X^{\ot 3}\ot P^{\ot 3})\geq
h^0(\omega _X^{\ot 2}\ot P^{\ot 2} \ot \pi ^*L)=$$
$$h^0(\pi _*(\omega _X^{\ot 2}\ot  P^{\ot 2})\ot L)
\geq \deg (\pi _*(\omega _X^{\ot 2} \ot P^{\ot 2}))+
rank(\pi _*(\omega _X^{\ot 2} \ot P^{\ot 2})) \ge 3+2$$
and this is the required contradiction.
\end{proof}

\begin{prop} \label{Iitaka}
Let $X$ be a smooth projective variety with $P_3(X)=4$, $q(X)=\dim
(X)$, and $f:X\lra Y$ be a birational model of its Iitaka
fibration. Then $Y$ is birational to an abelian variety.
\end{prop}

\begin{proof}
Since $X,Y$ are of maximal Albanese dimension, $K_{X/Y}$ is effective.
If $h^0(\omega _Y \ot P)>0$, it follows that $h^0(\omega _X \ot f^* P)>0$
and so by Proposition \ref{gt}, $f^*P=\OO _X$. By Proposition
\ref{albanese}, the map $f^*:\Pic ^0(Y)\lra \Pic ^0(X)$ is injective
and hence $P=\OO _Y$. Therefore $V^0(\omega _Y)=\{\OO _Y\}$ and by 
Theorem \ref{TCH2}, one has $\kappa (Y)=0$ 
and hence $Y$ is birational to an abelian variety.
\end{proof}

We are now ready to describe the cohomological support loci of
varieties with $\kappa (X) \ge 2$ explicitly. Recall that
by Proposition
\ref{PCH2}, for all
$\eta \ne 0 \in \bar{G}$, $\delta_\eta \ge 1$.

\begin{thm} Let $X$ be a smooth projective variety with
$P_3(X)=4, q(X)=\dim (X)$ and $\kappa(X) \ge 2$. Then $\kappa (X)=
2$ and  $\bar{G}\cong (\Z _2)^s$ for some $s\geq 1$.
\end{thm}
\begin{proof}
 The proof consists of following claims.
\begin{claim}\label{cl1}
If $\kappa (X)\geq 2$ and $T\subset V^0(\omega _X)$ is a positive
dimensional component, then $T+T\subset \Pic ^0(Y)$, i.e. $\bar{G}
\cong (\Z _2)^s$.
\end{claim}

\begin{proof}[Proof of Claim \ref{cl1}]
It suffices to prove that $2\eta =0$ for $0\ne \eta \in \bar{G}$. Suppose
that $2\eta \ne 0$, we will find a contradiction.

We first consider the case that $\delta_\eta \ge 2$ and
$\delta_{-2\eta} \ge 2$. Then by (\ref{pluri}), $P_{2,2\eta} \ge
1+1+\delta_\eta -1 \ge 3$, and $P_3 \ge 3+1+\delta_{-2\eta}-1 \ge
5$ which is impossible.

We then consider the case that $\delta_\eta \ge 2$ and
$\delta_{-2\eta} =1$.  Again we have  $P_{2,2\eta} \ge 3$. We
consider the induced map $\pi : X \to E:=T_{-2\eta}^\vee$ and
the inclusion $\varphi: \pi^*L \to \omega_X \ot Q_{-2\eta}$
where $E$ is an elliptic curve and $L$ is an ample line bundle on $E$.
It follows that there is an inclusion $$\pi^*L \ot (\omega _X\ot
Q_\eta)^{\ot 2} \to \omega_X^{\ot 3} \ot Q_\eta^{\ot 2} \ot 
Q_{-2\eta}.$$ By
Lemma \ref{elliptic}, one has that $rank (\pi_* (\omega_X \ot Q_\eta)^
{\ot 2}) \ge
2$. By Proposition \ref{Pm}, Riemann-Roch and Lemma \ref{L1}
$$ P_3 (X)=h^0(\omega_X^{\ot 3} \ot Q_\eta^{\ot 2} \ot Q_{-2\eta})
\ge h^0(\pi^*L \ot (\omega _X\ot Q_\eta)^{\ot 2}) =$$
$$ h^0((\omega _X \ot
Q_\eta)^{\ot 2})+rank 
(\pi_* (\omega _X\ot Q_\eta)^{\ot 2})\ge P_{2,2\eta}+2 \ge
5,$$ which is impossible.

Lastly, we consider the case that $\delta_\eta =1$. 
There is an induced map $\pi: X \to
E:=T_{\eta}^\vee$ and an inclusion $ \pi^*L \to \omega_X \ot
Q_{\eta}$. Hence there is an inclusion $\varphi: \pi^*L^{\ot 2} \to
(\omega_X \ot Q_{\eta})^{\ot 2}$. By Lemma \ref{elliptic}, we have
$P_{2,2\eta} \ge 3$. We now proceed as in the previous cases.

Therefore, any element $\eta \in \bar{G}$ is of order 2 and hence
$\bar{G}\cong (\Z_2)^s$.
\end{proof}

\begin{claim}\label{cldim}
If there is a surjective map with connected fibers to an elliptic curve
$\pi :X \lra E$ and an inclusion $\pi ^* L\lra \omega _X \ot P$ for
an ample line bundle $L$ on $E$ and $P\in \Pic ^0(X)$
(in particular if $\delta_i=1$ for some $i \ne 0$
cf. Corollary \ref{C9}). 
Then $\kappa(X)=2$.
\end{claim}
\begin{proof}[Proof of Claim \ref{cldim}]
Since $K_X$ is effective, there is also an
inclusion $L\lra \pi _*(\omega _X ^{\ot 2}\ot P)$. By Lemma
\ref{elliptic} , one has $rank (\pi _*(\omega _X ^{\ot 2}\ot
P))\geq 2$, $h^0(\pi _*(\omega _X ^{\ot 2}\ot P))\geq 2$. Consider
the inclusion $$\pi _*(\omega _X ^{\ot 2}\ot P)\ot L \lra \pi
_*(\omega _X ^{\ot 3}\ot P^{\ot 2}).$$ Since $$P_3(X)=h^0(\pi
_*(\omega _X ^{\ot 3}\ot P^{\ot 2}))\geq h^0( \pi _*(\omega _X
^{\ot 2}\ot P)\ot L)\geq$$ $$\deg (\pi _*(\omega _X ^{\ot 2}\ot
P))+ rank (\pi _*(\omega _X ^{\ot 2}\ot P)),$$ it follows that
$$\deg (\pi _*(\omega _X ^{\ot 2}\ot P))= rank (\pi _*(\omega _X
^{\ot 2}\ot P))=2$$ and the above homomorphism of sheaves induces
an isomorphism on global sections and hence is an isomorphism of
sheaves (cf. Proposition \ref{inclusion}). Therefore,
$$P_3(F)=h^0(\omega _F ^{\ot 3}\ot P^{\ot 2}) =2.$$ By Theorem \ref{T1},
it follows that $\kappa (F)=1$ and by easy addition, one has that
$$\kappa (X)\leq \kappa (F)+\dim (E)=2.$$
\end{proof}

\begin{claim}\label{cc1}
For all $i \ne 0$,  $P_{1,i} =1$.
\end{claim}

\begin{proof}[Proof of the Claim \ref{cc1} ]
If $P_{1,i}  \ge 2$, then by (\ref{pluri}), $$4 \ge P_2 \ge
2P_{1,i}+ \delta_i-1.$$ It follows that $\delta_i=1$. Let
$E=T^\vee$ and $\pi :X\lra E$ be the induced morphism. One has an
inclusion $\pi ^*L\lra \omega _X \ot Q_i$. By Lemma \ref{Lel}, one
has $h^0(E,L)=h^0(\omega _X \ot Q_i) \ge 2$. Consider the
inclusion $\pi ^*L^{\ot 2}\lra \omega_X^{\ot 2} \ot Q_i^{\ot 2}$. By Lemma
\ref{elliptic}, one sees that $$P_3 \ge P_{2,2i}=h^0(\omega_X^{\ot 2}
\ot Q_i^{\ot 2})> h^0(E, L^{\ot 2}) \ge 4,$$ which is impossible. 
\end{proof}

\begin{claim}\label{cc2} If $\kappa (X)=\dimm (S)$
for some component $S$ of $V^0(\omega _X)$, 
then $\kappa(X)=2$.
\end{claim}

\begin{proof}[Proof of Claim \ref{cc2}] Let 
$Q$ be a general point in $S$,
and $T$ be the translate of $S$ through the origin.
By Proposition \ref{Iitaka}, one sees that the induced map $X \to
T^\vee$ is isomorphic to the Iitaka fibration. We therefore
identify $Y$ with $T^\vee$. We assume that $\dimm (S)\geq 3$
and derive a contradiction. First of all, by (\ref{pluri})
$$P_3(X)=h^0(\omega _X ^{\ot 3}\ot Q^{\ot 2})\geq 
h^0(\omega _X ^{\ot 2}\ot Q)+\dimm (S)$$
and so $h^0(\omega _X ^{\ot 2}\ot Q)=1$ and $\dimm (S)=3$.

Let $H$ be an ample line bundle on $Y$ and for $m$ a sufficiently
big and divisible integer, fix a divisor $B\in |mK_X-f ^* H|$.
After replacing $X$ by an appropriate birational model, we may
assume that $B$ has simple normal crossings support. Let $L=\omega
_X\ot \OO _X (-\lfloor B/m\rfloor )$, then $L\equiv f ^* (H/m)+\{
B/m\}$ i.e. $L$ is numerically equivalent to the sum of the pull
back of an ample divisor and a k.l.t. divisor and so one has
$$h^i(Y, f_* (\omega _X \ot L \ot Q)\ot \eta )=0\ \ \ \ for\  all\
i>0\ \  and \ \eta \in \Pic ^0(Y).$$ Comparing the base loci, one can
see that $h^0(\omega _X \ot L\ot Q) =h^0(\omega _X ^{\ot 2}\ot
Q)=1$ (cf. \cite{CH1} Lemma 2.1 and Proposition \ref{Pm}) and so
$$h^0(Y,f _* (\omega _X \ot L \ot Q)\ot \eta )= h^0(f _* (\omega _X
\ot L \ot Q))=1\ \ \ \forall \eta \in \Pic ^0(Y).$$  Since
$f_*(\omega_ X\ot L\ot Q )$ is a torsion free sheaf of generic
rank one, by \cite{Hac} it is a principal polarization $M$.

Since one may arrange that $\lfloor \frac{B}{m} \rfloor \prec
K_X$. There is an inclusion $ \omega_X \ot Q \hookrightarrow
\omega_ X\ot L\ot Q$. Pushing forward to $Y$, it induces an inclusion
$$ \varphi: f_*(\omega_X \ot Q) \hookrightarrow M.$$
Since $f_*(\omega_X \ot Q)$ is torsion free, it is generically of rank
one. Hence it is of the form $M \ot \III_Z$ for some ideal sheaf
$\III_Z$. However, $h^0(Y,f_*(\omega_X \ot Q) \ot P)=h^0(M \ot P\ot 
\III_Z)
>0$ for all $P \in \Pic ^0(Y)$ and $M$ is a principal
polarization. It follows that $\III_Z=\OO_Y$ and thus $
f_*(\omega_X \ot Q)=M$. Therefore, one has an inclusion
$$ f^*M^{\ot 2} \hookrightarrow (\omega_X \ot Q) \ot  (\omega_X \ot L \ot
Q ) \hookrightarrow \omega_X^{\ot 3} \ot Q^{\ot 2}.$$ It follows that
$$4=P_3(X) =h^0(X,\omega_X^{\ot 3} \ot Q^{\ot 2}) \ge h^0(Y,M^{\ot 2}) \ge
2^{\dimm (S)}.$$ This is the required contradiction.
\end{proof}

\begin{claim}\label{cc4}
Any two components of 
$V^0(\omega _X)$ of dimension at least $2$ must be parallel.
\end{claim}
\begin{proof}[Proof of Claim \ref{cc4}]

For $i=1,2$, let $S_i:=T_i^\vee$ and $p_i:X\lra S_i$ be the
induced morphism. Assume that $\delta _1,\delta _2\geq 2$
and $T_1,T_2$ are not parallel.
By Lemma \ref{Lel}, one may write $K_X=G_i+D_i$
where $D_i$ is vertical with respect to $p_i:X\lra S_i$ and for
general $P\in T_i$, one has $|K_X+P|=G_i+|D_i+P|$ is a
$0$-dimensional linear system (see Claim \ref{cc1}).

Recall that we may assume that the image of the Iitaka fibration
$f:X\lra Y$ is an abelian variety.
Pick $H$ an ample divisor on $Y$ and for
$m$ sufficiently big and divisible integer, let
$$B\in |mK_X-f^* H|.$$
After replacing $X$ by an appropriate birational model, we may
assume that $B$ has normal crossings support.
Let $$ L:= \omega _X(-\lfloor \frac {B}{m}\rfloor )\equiv
\{ \frac{B}{m}\} +f^*H.$$
It follows that $$h^i(f_{*}(\omega _X \ot L \ot P )\ot \eta )=0\ \
\ for \ all \ i>0,\ \eta \in \Pic ^0(Y),\ \ P\in \Pic ^0(X).$$ The
quantity $h^0(\omega _X \ot L \ot P \ot f^*\eta )$ is independent
of $\eta \in \Pic ^0(Y)$. For some fixed $P\in T_1$ as above, and
$\eta \in \Pic ^0(S_1)$, one has a morphism
$$|D_1+P+\eta |\times |D_1+P-\eta |\lra |2D_1+2P|$$
and hence $h^0(\OO _X (2D_1)\ot P^{\ot 2})\geq 3$. Similarly
for some fixed $Q\in T_2$, and $\eta '\in \Pic ^0(S_2)$, one has a
morphism
$$|D_2+Q+\eta '|\times |K_X+L-Q+2P-\eta ' |\lra |K_X+L+D_2+2P|$$
and hence $h^0(\omega  _X (D_2)\ot L\ot P^{\ot 2})\geq 3$. It
follows that since $h^0(\omega _X ^{\ot 3}\ot P^{\ot 2})=4$, there
is a 1 dimensional intersection between the images of the 2
morphisms above which are contained in the loci
$$|2D_1+2P|+2G_1+K_X,\ \ \ \ \ 
|K_X+L+D_2+2P|+\lfloor \frac {B}{m}\rfloor+G_2.$$
It is easy to see that for all but finitely many $P\in \Pic
^0(X)$, one has $h^0(\omega _X \ot P)\leq 1$. So there is a 1
parameter family $\tau _2\subset \Pic ^0(S_2)$ such that for $\eta
' \in \tau _2$, one has that the divisor $D_{Q+\eta '}=|D_2+Q+\eta
'|$ is contained in $D_{P+\eta }+D_{P-\eta}+2G_1+K_X$ where $\eta
\in \tau _1$ a 1 parameter family in $\Pic ^0(S_1)$. Let
$D^*_{Q+\eta '}$ be the components of $D_{Q+\eta '}$ which are not
fixed for general $\eta ' \in \tau _2$, then $D^*_{Q+\eta '}$ is
not contained in the fixed divisor $2G_1+K_X$ and hence is
contained in some divisor of the form $D^*_{P+\eta }+D^*_{P-\eta
}$ and hence is $S_1$ vertical. 

If $\Pic ^0(S_1)\cap \Pic
^0(S_2)=\{ \OO _X\}$, then $D^*_{Q+\eta '}$ is $\alb$-exceptional,
and this is impossible by Lemma \ref{L7}. 

If there is a
1-dimensional component $\Gamma \subset \Pic ^0 (S_1) \cap \Pic
^0(S_2)$. Let $E=\Gamma^\vee $ and $\pi :X\lra E$ be the induced
morphism. The divisors $D^*_{Q+\eta '}$ are $E$-vertical. We may
assume that $\pi$ has connected fibers. Since the $D^*_{Q+\eta '}$
vary with $\eta ' \in \tau _2$, for general $\eta '\in \tau _2$,
they contain a smooth fiber of $\pi$. So for general $\eta '\in
\tau _2$ there is an inclusion $\pi ^* M\lra \omega _X\ot Q\ot \pi
^*\eta '$ where $M$ is a line bundle of degree at least 1.
By Claim \ref{cldim}, one has $\kappa (X)=2$ and hence $T_1,T_2$ are parallel.

If there is a
2-dimensional component $\Gamma \subset \Pic ^0 (S_1) \cap \Pic
^0(S_2)$, then $\delta _1=\delta _2\geq 3$. 
By (\ref{pluri}), one sees that $P_{2,Q_1+Q_2}\geq 3$.
By Lemma \ref{P2}, this is impossible.
\end{proof}
By Claim \ref{cldim}, if there is a one dimensional component,
then $\kappa (X)=2$. Therefore, we may assume that $\delta _i\geq 2$ 
for all $i\ne 0$. By Claim \ref{cc4}, since $\delta _i\geq 
2$ for all $i\ne 0$, 
then $S_i,S_j$ are parallel for all $i,j\neq 0$. 
By Theorem \ref{TCH2}, for an appropriate $i\ne 0$, 
$\kappa (X)=\dimm (S_i)$ and so by Claim \ref{cc2},
one has $\kappa (X)=2$.
\end{proof}

\section{Varieties of $P_3(X)=4,q(X)=\dimm (X)$ and $\kappa(X)=2$}
In this section, we classify varieties with 
$P_3(X)=4,q(X)=\dimm (X)$ and $\kappa(X)=2$.
The first first step is to describe the cohomological support loci
of these varieties.
We must show that the only possible cases are
the following (which corresponds to Examples 2 and 3
respectively):
\begin{enumerate}
\item $\bar{G} \cong \Z_2, V_0(X, \omega_X)=\{\OO_X\} \cup
S_\eta$, $\delta_\eta=2$.

\item $\bar{G} \cong \Z_2^2, V_0(X, \omega_X)=\{\OO_X\} \cup
S_\eta \cup S_\zeta \cup S_{\eta+\zeta}$,
$\delta_\eta=\delta_\zeta=1,\delta_{\eta+\zeta}=2$.
\end{enumerate}
Using this information, we will determine the sheaves $\alb _*(\omega _X)$
and this will enable us to prove the following:

\begin{thm}\label{Tk2}
 Let $X$ be a smooth projective variety with
$P_3(X)=4$, $q(X)=\dim (X)$ and $\kappa (X)=2$, then $X$ is
one of the varieties described in Examples 2 and 3.
\end{thm}
\begin{proof}
Recall that $f:X\lra Y$ is a morphism birational to the
Iitaka fibration, $Y$ is an abelian surface and $f=q\circ \alb$ where
$q:\Alb\lra Y$.
\begin{claim}\label{cl5}
One has that $f_*\omega _X=\OO _Y$.
\end{claim}
\begin{proof}[Proof of Claim \ref{cl5}]
By Proposition \ref{gt}, one has that
$V^0(\omega _X)\cap f^*\Pic ^0(Y)=\{\OO _X\}$.
By the proof of \cite{CH3} Theorem 4, one sees that
$f_*\omega _X\cong \OO _Y
\ot H^0(\omega _X)$. Since $h^0(\omega _X|F_{X/Y})=1$, it follows that
$rank(f_*\omega _X)=1$ and hence $f_*\omega _X\cong \OO _Y$.
\end{proof}

\begin{claim}\label{cl6} Let $T_1,T_2$ be distinct components of
$V^0(\omega _X)$ such that $T_1\cap T_2\ne \emptyset$, then
$T_1\cap T_2=P$ and $$f_*(\omega _X \ot P)=L_1\boxtimes L_2 \ot
\III _p$$ where $Y=E_1\times E_2$ and $L_i$ are line bundles of
degree $1$ on the elliptic curves $E_i$ and $p$ is a point of $Y$.
\end{claim} \begin{proof} [Proof of Claim \ref{cl6}]
Assume that $P\in T_1\cap T_2$. Since
$\kappa (X)=2$, by Proposition \ref{PCH2},
the $T_i$ are $1$-dimensional. Let $\pi _i:X\lra E_i:=T_i^\vee$ be
the induced morphisms. There are line bundles of degree 1, $L_i$
on $E_i$ and inclusions $\pi _i^*L_i\lra \omega _X \ot P$ (cf.
Corollary \ref{C9}).

We claim that $rank (\pi _{1,*}(\omega _X \ot P))=1$ . If this were not 
the case, then by Lemma \ref{fiber}
$$P_1(F_{X/E_1})=rank (\pi _{1,*}(\omega _X \ot P))\geq 2, 
\ \ \ \ P_2(F_{X/E_1})=rank (\pi _{1,*}(\omega _X^{\ot 2} \ot P))\geq 3$$
and so $$P_3(X)=h^0(\omega _X^{\ot 3}\ot P^{\ot 2})\geq h^0(\omega
_X^{\ot 2}\ot P\ot \pi _1^*L_1)=$$ $$h^0(\pi _{1,*}(\omega _X^{\ot
2}\ot P)\ot L_1)\geq rank(\pi _{1,*}(\omega _X^{\ot 2}\ot P))+
\deg (\pi _{1,*}(\omega _X^{\ot 2}\ot P))$$ and therefore
$$rank(\pi _{1,*}(\omega _X^{\ot 2}\ot P))=3,\ \ \ \ \ \deg (\pi
_{1,*}(\omega _X^{\ot 2}\ot P))=1.$$ Since $rank (\pi
_{1,*}(\omega _X ))=rank (\pi _{1,*}(\omega _X \ot P))$, one has
$$ \deg (\pi _{1,*}(\omega _X^{\ot 2}\ot P))\geq
 \deg ( \pi _{1,*}(       \omega _X)\ot L_1)\geq rank (\pi _{1,*}(\omega _X ))
\geq 2,$$
which is impossible.
Therefore, we may assume that $$rank (\pi _{i,*}(\omega _X \ot P))=1\ \
\ \ for\ i=1,2.$$
For any $P_i\in T_i$, one has that $P_i\ot P^\vee =\pi _i^* \eta _i$
with $\eta _i\in \Pic ^0(E_i)$. One sees that
$$h^0(\omega _X \ot P_i)=h^0(\pi _{i,*}(\omega _X \ot P)\ot \eta _i)=
h^0(\pi _{i,*}(\omega _X \ot P))=h^0(\omega _X \ot P).$$
If $h^0(\omega _X \ot P)\geq 2$, then we may assume that
$L_1:=\pi _{1,*}(\omega _X \ot P)$
is an ample line bundle of degree at least 2. From the inclusion
$\phi :L_1^{\ot 2}\lra \pi _{1,*}(\omega _X^{\ot 2}\ot P^{\ot 2})$,
one sees that $h^0(\omega _X^{\ot 2}\ot P^{\ot 2})=4$ and $\phi$ is an I.T. 0
isomorphism (cf. Lemma \ref{L1}) and so
$$P_2(F_{X/E_1})=h^0(\omega _X^{\ot 2}\ot P^{\ot 2}|F)=1.$$
By Theorem \ref{T1}, $\kappa (F_{X/E_1})=0$ and hence by easy addition,
$\kappa (X)\leq 1$ which is impossible. Therefore we may assume that
$h^0(\omega _X \ot P)=1$.

The coherent sheaf $f_*(\omega _X \ot P)$ is torsion free of generic rank 1
on $Y$ and hence is isomorphic to $L\ot \III$ where $L$ is a line bundle
and $\III$ is an ideal sheaf cosupported at finitely many points.
Let $q_i :Y\lra E_i$, so that $\pi _i = q_i \circ f$.
Since $$1=rank (\pi _{i,*}(\omega _X \ot P))=rank
(q_{i,*}(L\ot \III ) )=rank (q_{i,*}L),$$
one sees that $L.F_{Y/E_i}=1$ and it easily follows that
$L=L_1\boxtimes L_2$ where $L_i=q_{i,*}(L)$ is
a line bundle of degree 1 on $E_i$.
Clearly, $\III$ is the ideal sheaf of a point.
\end{proof}

We will now consider the case in which $\bar{G}=\Z _2$.
Let $B$ be the branch locus of $\alb : X\lra \Alb$.
The divisor $B$ is vertical with respect to $q:\Alb
\lra Y$ and hence we may write
$B=q^* \bar{B}$. Let $g\circ h:X\lra Z\lra A$ be the Stein factorization
of $\alb$. Then
$Z$ is a normal variety and $g$ is finite of degree $2$ and so
$g_*\OO _Z =\OO _A \oplus M^\vee$ where $M$ is a line bundle and
the branch locus $B$ is a divisor in $|2M|$.
The map $F_{Z/Y}\lra F_{\Alb /Y}$ is \'etale of degree 2 and so $M=q^*L
\ot P$ where $P$ is a $2$-torsion element of $\Pic ^0(X)$.
Let $\nu :\Alb '\lra \Alb$
be a birational morphism so that $\nu ^* B$ is a divisor with simple
normal crossings support. Let $B'=\nu ^*B-2\lfloor \nu ^*B/2\rfloor$
and $M'=\nu ^* (M)(-\lfloor \nu ^*B/2\rfloor)$.
Let $Z'$ be the normalization of $Z\times _\Alb \Alb '$, and
$g':Z'\lra \Alb '$
be the induced morphism. Then $g'$ is finite of degree 2,
$Z'$ is normal with rational singularities and
$g^\prime _*(\OO _{Z'})=\OO _{\Alb '}\oplus (M')^\vee$.
Let $\tilde {X}$ be an appropriate birational model
of $X$ such that there are morphisms $\alpha: \tilde {X}\lra \Alb '$, $v:
\tilde {X}\lra X$, $\tilde {\alb }:\tilde {X}\lra \Alb$ and
$\beta :\tilde {X}\lra Z'$. For all $n\geq 0$,
one has that $\beta _*(\omega _{\tilde{X}}^
{\ot n})\cong \omega _{Z^\prime}^{\ot n}$.
It follows that
$$\alpha  _*(\omega _{\tilde{X}}^{\ot m} )=\omega _{\Alb ^\prime}^{\ot m}\ot
({M'}^{\ot m-1}\oplus {M'}^{\ot m}).$$
Therefore
$$\alb _* (\omega _X)=\tilde {\alb }_* (\omega _{\tilde {X}})=
\nu _*(\omega _{\Alb ^\prime}\oplus \omega _{\Alb ^\prime}\ot
M')=$$ $$\OO _\Alb \oplus \nu _*(\omega _{\Alb ^\prime}\ot \nu ^* (q^*L)
(-\lfloor \nu ^*\frac{B}{2}\rfloor))=\OO _\Alb \oplus
q^*L\ot P\ot  \III (\frac{B}{2}).$$

\begin{claim}\label{cl9}
If $\bar{G}=\Z_2$, then for any $P\in V^0(\omega _X)$, one has
$$f_*(\omega _X \ot P)\ne
L_1\boxtimes L_2 \ot \III _p$$
where $Y=E_1\times E_2$ and $L_i$ are ample line bundles of
degree 1 on $E_i$ and $p$ is a point of $Y$.
\end{claim}
\begin{proof}[Proof of Claim \ref{cl9}]
If $f_*(\omega _X \ot P)=
L_1\boxtimes L_2 \ot \III _p$, then ${B}/{2}$ is not log terminal.
By \cite{Hac3} Theorem 1, one sees that
since ${B}/{2}$ is not log terminal,
one has that $\lfloor  {B}/{2}
\rfloor \ne 0$ and this is impossible as then $Z$ is not normal.
\end{proof}

Combining Claim \ref{cl6} and Claim \ref{cl9}, one sees that if
$\bar{G}=\Z_2$, then $V_0(X,\omega_X)=\{\OO _X\} \cup S_\eta$ with
$\delta_\eta=2$. We then have the following:

\begin{claim} \label{cl10}
If $\bar{G}=\Z_2$, then $h^0(X,\omega_X \ot P)=1$ for all $P \in
S_\eta$.
\end{claim}
\begin{proof}[Proof of Claim \ref{cl10}]
It is clear that $h^0(\tilde{X},\omega _{ \tilde{X}}\ot P )= 
h^0(A',\omega_{A'}
\ot M' \ot  P)$ for all $P \in S_\eta$, and 
$h^0(\tilde{X},\omega _{ \tilde{X}}\ot P )=1$ for general
$P \in S_\eta$.

If $h^0(\tilde{X},\omega _{ \tilde{X}} \ot Q_0 ) \ge
2$ for some $Q_0\in S_\eta$, then $h^0(\tilde{X},
\omega  _{ \tilde{X}}\ot Q_0 )=2$ as otherwise
$h^0(\omega  _{ \tilde{X}}^{\ot 2}
\ot Q_0 ^{\ot 2}) \ge 3+3-1$ which is impossible.

Consider the linear series $|K_{A'}+M'+Q_0|$. Let $\mu:\tilde{A} \to A'$ be a
log resolution of this linear series. We have
$$ \mu^* |K_{A'}+M'+Q_0| = |D|+F,$$
where $|D|$ is base point free and $F$ has simple normal crossings support.
There is an induced map $\phi_{|D|}: \tilde{A} \to \PPP^1$
such that $|D|=\phi_{|D|}^*|\OO _{\PPP^1}(1)|$.
We have an inclusion
$$ \varphi_1: \phi_{|D|}^*|\OO _{\PPP^1}(2)|+G \hookrightarrow
\mu^*|2K_{A'}+2M'+2Q_0|.$$ For all $\eta \in \Pic ^0(Y)$, there is a
morphism
$$ \varphi_2: \mu^*|K_{A'}+M'+Q_0+\eta|+\mu^*|K_{A'}+M'+Q_0-\eta|
\lra \mu^*|2K_{A'}+2M'+2Q_0|.$$ Notice that
$h^0(A',\omega_{A'}^{\ot 2}\ot {M'}^{\ot 2}\ot Q^{\ot 2}_0) 
\le h^0(X,  \omega_X^{\ot 2}\ot Q^{\ot 2}_0)\le 4$.
Since $h^0(\PPP^1, \OO_{\PPP^1}(2))=3$, $\varphi _1$ has a
$2$-dimensional image. Since $\eta$ varies in a $2$-dimensional
family, $\varphi_2$ also has  $2$-dimensional image. In
particular, there is a positive dimensional family
$\mathcal{N}\subset \Pic ^0(Y)$ such that for general $\eta \in
\mathcal{N}$, one has $$D_{\pm
\eta} +F_{\pm \eta}\in \mu^*|K_{A'}+M'+Q_0 \pm \eta|$$ 
where $G=F_\eta +F_{-\eta}$ and $D_\eta
+D_{-\eta}\in \phi_{|D|}^*|\OO _{\PPP^1}(2)|$. Since $G$ is a
fixed divisor, it decomposes in at most finitely many ways as the
sum of two effective divisors and so we may assume that
$F_{\eta},F_{-\eta}$ do not depend on $\eta\in \mathcal{N}$.

Take any $\eta \ne \eta' \in \mathcal{N}$ with $F_\eta=F_{\eta'}$.
One has that $D_\eta=\phi_{|D|}^* H$ is numerically equivalent to
$D_{\eta'}=\phi_{|D|}^* H'$. It follows that $H$ and $H'$ are
numerically equivalent on $\PPP^1$ hence linearly equivalent. Thus
$D_\eta$ and $D_{\eta'}$ are linearly equivalent which is a
contradiction.

\end{proof}

\begin{claim}\label{cl11}
If $\bar{G}=\Z _2$, then $\alb :X\lra \Alb$ has generic degree 2
and is branched over a divisor $B\in |2f^*\Theta |$ where $\OO
_Y(\Theta)$ is an ample line bundle of degree 1. Furthermore,
$\alb _* (\OO _X)\cong \OO _\Alb \oplus q^* \OO _Y(\Theta )\ot P$
where $P\notin \Pic ^0(Y)$ and $P^{\ot 2}=\OO _\Alb$. See Example
2.
\end{claim}
\begin{proof}[Proof of Claim \ref{cl11}]

For all $\eta \in \Pic ^0(Y)$ and $P\in S_\eta$, 
one has that $$h^0(\omega _X \ot
P\ot \eta ) =h^0(\omega _{\Alb ^\prime }\ot M' \ot P \ot \eta
)=1.$$ 
The sheaf $q_* \nu _*(\omega _{\Alb '}\ot M'\ot P)$
is torsion free of generic rank $1$ and $$h^0(q_* \nu _*(\omega
_{\Alb '}\ot M'\ot P)\ot \eta )=1\ \ \  for\ all\ \eta \in \Pic ^0(Y).$$
Following the proof of Proposition 4.2 of \cite{HP}, one sees that
higher cohomologies vanish. By \cite{Hac}, $q_* \nu
_*(\omega _{\Alb '}\ot M'\ot P)$ is a principal polarization $\OO
_Y(\Theta )$. From the isomorphism $\nu _*(\omega _{\Alb '}\ot
M'\ot P)\cong \bar{L}\ot \III (\bar{B}/2)$, one sees that
$\bar{L}=\OO _Y (\Theta )$ and $\III (\bar{B}/2)=\OO _Y .$
Therefore, $ \nu _*(\omega _{\Alb '}\ot M'\ot P )\cong q^* \OO _Y
(\Theta )$. It follows that $$\alb _*(\omega _X)\cong \OO _\Alb \oplus q^*
\OO _Y(\Theta )\ot P.$$
\end{proof}

From now on we therefore assume that $\bar{G} \ne \Z_2$.

\begin{claim}\label{cl7} $V^0(K_X)$
has at most one $2$-dimensional component.
\end{claim}

\begin{proof} [Proof of Claim \ref{cl7}]
Let $S_\eta,S_\zeta$ be  $2$-dimensional components
of $V^0(\omega _X)$ with $\eta \ne \zeta$. Since $\kappa(X)=2$,
one has $\delta_{\eta,\zeta}=2$. Thus  by (\ref{pluri}),
$P_{2,\eta+\zeta} \ge 3$. By Lemma \ref{P2}, this is impossible.
\end{proof}

\begin{claim}\label{cl8}
Let $T_1,T_2$ be two parallel
$1$-dimensional components of $V^0(\omega _X)$, then $T_1+\Pic
^0(Y)=T_2+\Pic ^0(Y)$.
\end{claim}
\begin{proof}[Proof of Claim \ref{cl8}]

Let $P_i\in T_i$, $\pi :X\lra E:=T_1^\vee=T_2^\vee$ the induced
morphism and $L_i$ ample line bundles on $E_i$ with inclusions
$\phi _i: \pi ^* L_i\lra \omega _X \ot P_i$. By Lemma
\ref{claimA}, one sees that $h^0(\omega _X ^{\ot 2}\ot P_1\ot
P_2)\geq 2$. If it were equal, then the inclusion $$L_1\ot L_2
\lra \pi _* (\omega _X ^{\ot 2}\ot P_1\ot P_2)$$ would be an I.T.
0 isomorphisms and this would imply that $P_2(F_{X/E})=1$ and
hence that $\kappa (X)\leq 1$. So $h^0(\omega _X ^{\ot 2}\ot
P_1\ot P_2)\geq 3$. 
By Lemma \ref{P2}, this is impossible.
\end{proof}

\begin{claim}\label{cl12}
If $\bar {G}\ne \Z_2$, let $S_\eta$ be a $2$-dimensional component of
$V^0(\omega _X )$, then $h^0(\omega _X \ot P)=1$ for all $P\in S_\eta$.
In particular $f_*(\omega _X \ot P)$ is a principal polarization.
\end{claim}

\begin{proof}[Proof of Claim \ref{cl12}]
Let $f:X\lra (S_\eta )^\vee$ be the induced morphism.
Then $f$ is birational to the Iitaka fibration of $X$
i.e. $(S_\eta )^\vee =Y$.
By Claim \ref{cl7}, $V^0(\omega _X)$ has at most one $2$-dimensional
component, and so there must exist a $1$-dimensional component $S_\zeta$
of $V^0(\omega _X)$. Let $ \pi :X\lra E:=T_\zeta^\vee$
be the induced morphism. There is an ample line bundle $L$
on $E$ and an inclusion $\pi ^*L \lra \omega _X \ot Q_\zeta$ for
some general $Q_\zeta \in S_\zeta$.

Assume that $P\in S_\eta$ and
$h^0(\omega _X \ot P)\geq 2$.
If $rank (\pi _*(\omega _X \ot P))=1$, then $\pi _*(\omega _X \ot P)$
is an ample line bundle of degree at least $2$ and hence
$h^0( \pi _*(\omega _X \ot P)\ot \eta )\geq 2$ for all $\eta
\in \Pic ^0(E)$.
It follows that
$$h^0(\omega _X ^{\ot 2}\ot P\ot Q_\zeta)\geq h^0(\omega _X \ot P\ot \pi ^*L)=
h^0(\pi _*(\omega _X \ot P)\ot L)\geq 3.$$
By Lemma \ref{P2}, this is impossible.
Therefore, we may assume that 
$rank (\pi _*(\omega _X \ot P))\ge 2$. Proceeding as above , since 
$$h^0(\pi _*(\omega _X \ot P)\ot L)\geq
rank (\pi _* (\omega _X \ot P))+\deg  (\pi _* (\omega _X \ot P)),$$
it follows that 
$\pi _*(\omega _X\ot P)$ 
is a sheaf of degree $0$. Since $h^0(\pi _*(\omega _X\ot P)\ot \eta )
>0$ for all $\eta \in \Pic ^0(E)$, By Riemann-Roch one sees that
also $h^1(\pi _*(\omega _X\ot P)\ot \eta )>0$ for all $\eta \in \Pic ^0(E)$.
By Theroem \ref{genvanish},
this is impossible.

Finally, the sheaf $f_*(\omega _X \ot P)$ is torsion free of generic
rank $1$ on $Y$ and hence, by \cite{Hac}, it is a principal polarization.
\end{proof}

\begin{claim}\label{cl13}
Assume that $\bar {G} \ne \Z _2$. Then, for any $P\in V^0(\omega
_X)-\Pic ^0(Y)$ one has that $f_* (\omega _X\ot P)$ is either:

i) a principal polarization
on $Y$,

ii) the pull-back of a line bundle of degree 1 on an
elliptic curve or

iii) of the form $L\boxtimes L'\ot \III _p$ where $L,L'$
are ample line bundles of degree $1$ on $E,E'$, $Y=E\times E'$ and $p$
is a point of $Y$.

In particular, there are no 2 distinct
parallel components of $V^0(\omega _X)$.
\end{claim}
\begin{proof}[Proof of Claim \ref{cl13}]
By Claim \ref{cl12}, we only need to consider the case in which
all the components of
$(P+\Pic ^0(Y))\cap V^0(\omega _X)$ are $1$-dimensional.
By Claim \ref{cl6}, we may also assume that these components are parallel.

For any 1 dimensional component
$T_i$ of $(P+\Pic ^0(Y))\cap V^0(\omega _X)$,
$P_i\in T_i$ and corresponding projection
$\pi _i:X\lra E_i:=T_i^\vee$, one has $rank (\pi _{i,*}
(\omega _X \ot P_i))=1$ and hence $\pi _{i,*}
(\omega _X \ot P_i)={L_i}$ is an ample line bundle of degree
at least 1 on $E_i$. If this were not the case, then
By Lemma \ref{fiber}, $$rank (\pi _{i,*}
(\omega _X \ot P_i))=h^0(\omega _F)\geq 2$$ and so $$rank (\pi _{i,*}
(\omega _X ^{\ot 2}\ot P_i))=h^0(\omega _F^{\ot 2})\geq 3.$$
From the inclusion (cf. Corollary
\ref{C9})
$$\pi _i^*L_i\lra \omega _X \ot P_i\lra \omega _X^{\ot 2}\ot P_i,$$
one sees that $h^0(\omega _X^{\ot 2}\ot P_i)\geq 2$
(cf. Lemma \ref{elliptic}). By Lemma \ref{L1}, $\deg (
\pi _{i,*}(\omega _X^{\ot 2}\ot P_i))\geq 2$.
By Riemann-Roch, one has $$h^0(L\ot \pi _{i,*}(\omega _X^{\ot 2}\ot P_i))
\geq \deg (\pi _{i,*}(\omega _X^{\ot 2}\ot P_i))+rank
(\pi _{i,*}(\omega _X^{\ot 2}\ot P_i))\geq 5.$$
This is a contradiction and so 
$rank (\pi _{i,*}
(\omega _X \ot P_i))=1$.

Since we assumed that all components of $V^0(\omega _X)\cap (P+\Pic ^0(Y))$
are parallel, then one has $\pi _i=\pi$, $E=E_i$ are independent of $i$.
Let $q:Y\lra E$. Since there are injections 
$$\Pic ^0(E)+P_1=T_1\hookrightarrow P_1+\Pic ^0(Y)
\hookrightarrow \Pic ^0(X),$$  
we may assume that $q$ has connected fibers.
The sheaf $f_*(\omega _X \ot P_1)$
is torsion free of rank 1, and hence we may write
$f_*(\omega _X \ot P_1)\cong M\ot \III$ where $M$ is a line bundle
and $\III$ is supported in codimension at least 2 (i.e. on points).
Since $ rank (\pi _{*}
(\omega _X \ot P_1))=1$, one has that $h^0(M|F_{Y/E})=1$.

For general $\eta \in \Pic ^0(Y)$, one has that $V^0(\omega _X)
\cap P_1+\eta + \Pic ^0(E)=\emptyset$ and so the semipositive 
torsion free sheaf $\pi _*(\omega _X \ot P_1 \ot \eta )$
must be the $0$-sheaf. In particular $h^0(M\ot \eta |F_{Y/E})=0$.
It follows that $\deg (M|F_{Y/E})=0$ and hence $M|F_{Y/E}=\OO _{F_{Y/E}}$.
One easily sees that $h^0(M\ot \eta )=0$ for all $\eta \in \Pic ^0(Y)-
\Pic ^0(E)$ and hence $$V^0(\omega _X)=P_1+\Pic ^0(E)=T_1.$$
By Proposition \ref{inclusion}, one has that  $q ^* L_1$ and 
$ f_*(\omega _X \ot P_1)$ are isomorphic if and only if the inclusion
$q^*L_1\lra f_*(\omega _X \ot P_1)$ induces isomorphisms
$$H^i(Y,q ^* L_1 \ot \eta )\lra H^i(Y, f_*(\omega _X \ot P_1)\ot \eta)$$ 
for $i=0,1,2$ and all $\eta \in \Pic ^0(Y)$.
If $\eta \in \Pic ^0(Y)-\Pic ^0(E)$,
then both groups vanish and 
so the isomorphism follows. If $\eta \in \Pic ^0(E)$, we proceed as follows:
Let $p:\Alb \lra E$ and $W\subset H^1(A,\OO _A)$ a linear subspace 
complementary to the tangent space to $T_1$. By Proposition 2.12 of 
\cite{Hac2}, one has isomorphisms
$$H^1(\alb _*(\omega _X \ot P_1)\ot p^* \eta )\cong
H^0(\alb _*(\omega _X \ot P_1)\ot p^* \eta )\ot \wedge ^iW\cong$$
$$H^0(q^*(L_1\ot \eta))\ot \wedge ^iW\cong H^i(q^*L_1\ot \eta)).$$
Pushing forward to $Y$, one obtains the required isomorphisms.
\end{proof}

\begin{claim} \label{cl16}
If $\bar {G}\ne \Z _2$, then $\bar {G}= (\Z _2)^2$ and 
$$V_0(X,\omega_X)=\{\OO_X\} \cup
S_\eta \cup S_\zeta \cup S_\xi$$ with $\delta_\eta=2$,
$\delta_\zeta=\delta_\xi=1$.
\end{claim}
\begin{proof}[Proof of Claim \ref{cl16}]
We have seen that $V^0(\omega _X)$ has at most one
$2$-dimensional component and there are no parallel
$1$-dimensional components.
Since $\bar{G} \ne \Z_2$, then there are at least two $1$-dimensional
components of $V^0(\omega _X)$.
We will show that given two one dimensional components contained
in $Q_1+\Pic ^0(Y)\ne Q_2+\Pic ^0(Y)$, then $$\left( Q_1+Q_2+\Pic ^0(Y)
\right)
\cap V^0(\omega _X)$$ does not contain a 1-dimensional component.
Grant this for the time being.
Then, by Proposition \ref{PCH2},
it follows that $Q_1+Q_2+\Pic ^0(Y)$ is a $2$-dimensional component of
$V^0(\omega _X)$. If $|\bar {G}|>4$, this implies that there
are at least two $2$-dimensional components, which is impossible,
and so $|\bar{G}|=4$ and the claim follows.

Suppose now that there are three $1$-dimensional components of
$V^0(\omega _X)$, say $S_1,S_2,S_3$, contained in $Q_1+\Pic
^0(Y),Q_2+\Pic ^0(Y), Q_3+\Pic ^0(Y)$ respectively with 
$Q_1+Q_2+Q_3\in \Pic ^0(Y)$. By Claim \ref{cl13},
these components are not parallel to each other. We may
assume that $\pi_i: X \to E_i:=S_i^\vee$
 factors through $f: X \to Y$ and that $Y$ is an abelian surface.
Let $q_i:Y\lra E_i$ be the induced morphisms.

Let $Q_1,Q_2,Q_3$ be general torsion elements in $S_1,S_2,S_3$ and
$$\GG:=f_* (\omega _X^{\ot 2}\ot Q_2\ot Q_3),\ \ \ \FF:=
f_* (\omega _X^{\ot 3}\ot Q_1\ot Q_2\ot Q_3 ).$$
From the inclusions $\pi _i^*L_i\lra \omega _X \ot Q_i$,
one sees that we have
inclusions $$\varphi: q_2^*L_2 \otimes q_3^*L_3 \to \GG,\ \ \
\psi :q_1^*L_1 \otimes q_2^*L_2\otimes q_3^*L_3 \to \FF$$ where $L_i$ are
ample line bundles on $E_i$ respectively.
Since $\FF$ is torsion free of generic rank one, we may write
$$\FF = q_1^*L_1 \otimes q_2^*L_2\otimes q_3^*L_3\ot N \ot \III $$
where $N$ is a semi-positive line bundle on $Y$ and $\III$
is an ideal sheaf cosupported at points.
If $N$ is not numerically trivial (or if $F_{Y/E_1}\cdot 
q_i^*L_i>1$ for $i=2$ or $i=3$), then $N$ is not vertical with respect
to one of the projections $q_i$, say $q_1$. Then
$$rank (q_{1,*}(\FF ))= F_{Y/E_1}\cdot ( q_1^*L_1 + q_2^*L_2+ q_3^*L_3+ N)
\geq 3.$$
On the other hand, from the inclusion $\varphi$, one sees that
$rank (q_{1,*}(\GG ))\geq 2$. Consider the inclusion 
of I.T. 0 sheaves $L_1\lra
q_{1,*}(\GG \ot \eta )$ with $\eta =Q_1\ot Q_2^\vee \ot Q_3^\vee\in
\Pic ^0(Y)$. Since it is not an isomorphism,
one sees that $$h^0(\GG )=h^0(\GG \ot \eta) > h^0(L_1)\geq 1.$$
From the inclusion
$$\rho : L_1\ot q_{1,*}(\GG )\lra q_{1,*}(\FF )=\pi _{1,*}
(\omega _X^{\ot 3}\ot Q_1\ot Q_2\ot Q_3 )$$
one sees that by Riemann-Roch
$$h^0(\GG )+rank (q_{1,*}(\GG ))\leq 
h^0(\omega _X^{\ot 3}\ot Q_1\ot Q_2\ot Q_3 )=P_3(X)$$
and therefore
$$h^0(\GG )=2,\ \ \ rank (q_{1,*}(\GG ))=2.$$
In particular, $\rho$ is an I.T. 0 isomorphism. So,
$rank (q_{1,*}(\FF )) =rank  (q_{1,*}(\GG ))=2$ which is a contradiction.
Therefore, we have that $$N\in \Pic ^0(Y)\ \ \ and
\ \ \ q_2^*L_2.F_{Y/E_1}=q_3^*L_3.F_{Y/E_1}=1.$$ Since $deg (L_i)=1$,
one has $q_i^*L_i\equiv
F_{Y/E_i}$. Since
$(q_1^*L_1 \otimes q_2^*L_2\otimes q_3^*L_3)^{ 2}\geq 8$,
we have that $q_2^*L_2\cdot q_3^*L_3\geq 2$.
Since
$$h^0(q_2^*L_2\otimes q_3^*L_3)\leq h^0(\GG )=2,$$
one sees that $q_2^*L_2.q_3^*L_3=2$ and hence $\III =\OO _Y$.

Now let $\GG ':=f_* (\omega _X^{\ot 2}\ot Q_1\ot Q_3)$.
Proceeding as above, one sees that $$rank (q_{2,*}\GG ')\geq
F_{Y/E_2}\cdot (q_1^*L_1 +q_2^*L_3)=3,\ \ \ \ h^0(q_{2,*}\GG ')>
h^0(L_2)=1.$$
By Riemann Roch, one has that
$$P_3(X)=h^0(\omega _X ^{\ot 3}\ot Q_1\ot Q_2 \ot Q_3 )\geq
h^0(L_1 \ot q_{2,*}\GG ')\geq 5$$
which is the required contradiction.
\end{proof}

\begin{claim}\label{cl15}
If $\bar {G}\cong (\Z _2)^2$, then $Y=E_1\times E_2$ and there are
line bundles $L_i$ of degree $1$ on $E_i$, projections $p_i: \Alb
\lra E_i$ and $2-$torsion elements $Q_1,Q_2\in \Pic ^0(X)$ that
generate $\bar {G}$, such that
$$\alb _* (\OO _X)\cong \OO _\Alb \oplus M_1^\vee
\oplus M_2^\vee \oplus M_1^\vee \ot M_2^\vee  $$
with
$$M_1=p _1^*L_1\ot Q_1^\vee ,\ \ \ M_2=p _2^* L_2 \ot Q _2^\vee \ \ and\ \ \
M_3= M_1 \ot M_2.$$
In particular $X$ is birational to the fiber product of two
degree 2 coverings $X_i\lra \Alb$ with $P_3(X_i)=2$.
\end{claim}
\begin{proof}[Proof of Claim \ref{cl15}]
By Claim \ref{cl16}, the degree of $\alb :X\lra \Alb$ is
$|\bar{G}|=4$ and there are two
non parallel $1$-dimensional
components of $V^0(\omega _X)$ say $S_1,S_2$ such that $S_1+\Pic ^0(Y)\ne
S_2+\Pic ^0(Y)$. Let $E_i:=S_i^\vee$ and
$q _i:Y\lra E_i$, $\pi _i :X\lra E_i$
be the induced morphisms. Then there are inclusions
$\pi_i^*L_i\lra \omega _X \ot Q_i$ where $Q_i\in S_i$.
Moreover, by Claim \ref{cl16},
$Q_1+Q_2+\Pic ^0(Y)\subset V^0(\omega _X)$. By Claim \ref{cl12}, one has that
$$L:=f_*(\omega _X \ot Q_1 \ot Q_2)$$
is an ample line bundle of degree $1$. 
Moreover,
$$V^0(\omega _X)=\{\OO _X\}\cup S_1 \cup S_2\cup (Q_1+Q_2+\Pic ^0(Y)).$$
From the inclusion $$q_1^*L_1\ot q_2^*L_2\ot L\lra
f_*(\omega _X ^{\ot 3}\ot Q_1^{\ot 2}\ot Q_2^{\ot 2})$$
and the equality $4=P_3(X)=h^0(\omega _X ^{\ot 3}\ot Q_1^{\ot 2}\ot
Q_2^{\ot 2})$, one sees that $$L^2=2,\ \ \  L .q_i^*L_i=q_1^*L_1.q_2^*
L_2=1.$$ By the Hodge Index Theorem, one sees that since
$$L^2(q_1^*L_1+q_2^*L_2)^2=\left( L.(q_1^*L_1+q_2^*L_2)\right)^2$$
then the principal polarization $L$ is numerically 
equivalent to $q_1^*L_1+q_2^*L_2$. Therefore,
$$(Y,q_1^*L_1\ot q_2^*L_2)\cong (E_1,L_1)\times (E_2,L_2),$$
and one sees that 
$$L=q_1^*(L_1\ot P_1)\ot q_2^*(L_2\ot P_2),\ \ \ P_i\in \Pic ^0(E_i).$$
We have inclusions $$L\lra f_*(\omega _X \ot Q_1
\ot Q_2)\lra f_*(\omega _X ^{\ot 2}\ot Q_1
\ot Q_2),$$ $$q_1^*L_1\otimes q_2^*
L_2\lra f_*(\omega _X ^{\ot 2}\ot Q_1
\ot Q_2).$$ Let $\GG :=\omega _X ^{\ot 2}\ot Q_1
\ot Q_2$. If $h^0(\GG )=1$, then $L=q_1^*L_1\ot q_2^*L_2$ as required.
If $h^0(\GG )\geq 2$, then one sees that
$$h^0(\pi _{1,*}(\GG )\ot L_1\ot P_1)\geq rank (\GG )+\deg (\GG )\geq 1+2.$$
Since $$rank (\pi _{2,*}(\GG \ot \pi _1^*(L_1\ot P_1)))\geq
rank (q_{2,*}(q_1^*(L_1^{\ot 2}\ot P_1)\ot q_2^*(L_2)))=2,$$
one sees that
$$P_3(X)\geq h^0(\omega _X^{\ot 2}\ot Q_1\ot Q_2 \ot L)=
h^0(\pi _{2,*}(\GG \ot \pi _1^*(L_1\ot P_1)) \ot L_2 \ot P_2)\geq 2+3$$
and this is impossible.
Let $M_i:=p_i^*L_i\ot Q_i^\vee$. By Claim \ref{cl13}, one has
$$\alb _*(\omega _X)\cong \OO _\Alb \oplus M_1\oplus M_2 \oplus M_1 \ot M_2$$
and hence by Groethendieck duality,
$$\alb _* (\OO _X)\cong \OO _\Alb \oplus M_1^\vee
\oplus M_2^\vee \oplus M_1^\vee \ot M_2^\vee . $$
Let $X\lra Z\lra \Alb$ be the Stein factorization.
Following \cite{HM} \S 7, one sees that the only possible nonzero
structure constants defining the $4-1$ cover $Z\lra \Alb$ are
$c_{1,4}\in H^0(M_1\ot M_2\ot M_3^\vee)$, $c_{1,6}\in H^0(M_1\ot
M_2^\vee \ot M_3)$ and $c_{4,6}\in H^0(M_1^\vee \ot M_2 \ot M_3)$.
So, $Z\lra \Alb$ is a
bi-double cover. It is determined by two degree $2$ covers $\alb_i:X_i \lra
\Alb$
defined by $\alb_{i,*}(\OO _{X_i})=\OO _\Alb \oplus p_i^*L_i\ot Q_i^\vee $
and sections $-c_{1,4}c_{1,6}\in H^0(M_1^{\ot 2})$ and $c_{1,4}c_{4,6}
\in H^0(M_2^{\ot 2})$. It is easy to see that $X_1,X_2,Z$
are smooth.
\end{proof}

This completes the proof.
\end{proof}

\section{Varieties with $P_3(X)=4,q(X)=\dimm (X)$ and  $\kappa(X)=1$}
\begin{thm} \label{main} Let $X$ be a smooth projective
variety with $P_3(X)=4$, $q(X)=\dim (X)$ and $\kappa (X)=1$ then
$X$ is birational to $(C\times \tilde{K})/G$ where $G$ is an
abelian group acting faithfully by translations on an abelian variety
$\tilde{K}$ and faithfully
on a curve $C$. The Iitaka fibration of $X$ is birational to
$f:(C\times \tilde{K})/G\lra C/G=E$ where $E$ is an elliptic curve
and $\dim H^0(C,\omega _C^{\ot 3})^G=4$. \end{thm} \begin{proof} Let
$f:X\lra Y$ be the Iitaka fibration. Since $\kappa (X)=1$, 
and $\alb:X\lra \Alb$ is generically finite, one has
that $Y$ is a curve of genus $g\geq 1$. If $g=1$, then $Y$ is an
elliptic curve and by 
Proposition \ref{albanese},
$Y\lra \Alb (Y)$ is of degree 1 (i.e. an isomorphism). By
Proposition \ref{albanese} one sees that if $g\geq 2$, then
$q(X)\geq \dim (X)+1$ which is impossible.

From now on we will denote the elliptic curve $\Alb (Y)$ simply by $E$
and $f:X\lra E$ will be the corresponding algebraic fiber space.
Let $X\lra \bar{X}\lra \Alb$ be the Stein factorization of the Albanese map. 
Since
$\bar{X}\lra \Alb$ is isotrivial, there is a generically finite
cover $C \lra E$ such that $\bar{X}\times _E C$ is birational to
$C\times \tilde{K}$. We may assume that $C\lra E$ is a Galois
cover with group $G$.
$G$ acts by translations on $\tilde {K}$ and we may assume that
the action of $G$ is faithful on $C$ and $\tilde {K}$. Since $G$
acts freely on $C\times \tilde{K}$, one has that $$H^0(X, \omega
_X^{\ot 3})=H^0(C\times \tilde{K}, \omega _{C\times \tilde{K}}
^{\ot 3})^G=H^0(\tilde{K}, \omega_{\tilde{K}}^{\ot 3})  \otimes
H^0(C,\omega_C^{\ot 3})]^G.$$ Since $G$ acts on $\tilde{K}$ by
translations, $G$ acts on $H^0(\tilde{K}, \omega_{\tilde{K}}^{\ot
3})$ trivially. It follows that $$4=P_3(X)=\dim
H^0(C,\omega_C^{\ot 3})^G.$$ Similarly, one sees that
$q(X)=q(C/G)+q(\tilde{K}/G)$ and so $q(C/G)=1$.
\end{proof}

We now consider the induced morphism $\pi: C \to C/G=:E$. By the
argument of \cite{Be}, Example VI.12, one has $$4=\dim
H^0(C,\omega_C^{\ot 3})^G=h^0(E, \OO(\sum_{P \in E} \lfloor 3
(1-\frac{1}{e_P}) \rfloor)).$$ Where  $P$ is a branch points of
$\pi$, and $e_P$ is the ramification index of a ramification point
lying over $P$. Note that $|G|=e_P s_P$, where $s_P$ is the number
of ramification points lying over $P$.

It is easy to see that since $$\lfloor 3 (1-\frac{1}{e_P}) \rfloor
= 1\ (\text{resp.}\ =2) \ \text{ if }\  e_P=2 \ (\text{resp.}\ e_P
\ge 3),$$ we have the following cases:

\medskip \noindent {\bf Case 1.} $4$ branch points $P_1,..,P_4$
with $e_{P_i}=2$.\\ {\bf Case 2.} $3$ branch points $P_1,P_2,P_3$
with $e_{P_1}\ge 3, e_{P_2}=e_{P_3}=2$.\\ {\bf Case 3.} $2$ branch
points $P_1,P_2$ with $e_{P_i}\ge 3$.\\

We will follow the notation of \cite{Pa}. Let $\pi: C \to E$ be an
abelian cover with abelian Galois group $G$. There is a splitting
$$ \pi_* \OO_C = \oplus_{\chi \in G^*} L^{\vee} _\chi.$$ In
particular, if $d_\chi :=\deg (L_{\chi})$, then $$g=1+\sum_{\chi
\in G^*,\  \chi \ne 1} d_\chi.$$

For every branch point $P_i$ with $i=1,...,s$, the inertia group
$H_i$, which is defined as the stabilizer subgroup at any point
lying over $P_i$, is a cyclic subgroup of order $e_{i}:=e_{P_i}$.
We also associate a generator $\psi_i$ of each $H_i^*$ which
corresponds to the character of $P_i$. For every $\chi \in G^*$,
$\chi_{|H_i}=\psi_i^{n(\chi)}$ with $0 \le n(\chi) \le |H_i|-1$. And
define $$\epsilon^{H_i,\psi _i}_{\chi,\chi'}:= \lfloor
\frac{n(\chi)+n(\chi')}{|H_i|} \rfloor.$$ Following \cite{Pa}, one
sees that there is an abelian cover $C \to E$ with group $G$ with
building data $L_\chi$ if and only if the line bundles $L_\chi$
satisfy the following set of linear equivalences:
\begin{equation}
 L_\chi+L_{\chi'} = L_{\chi \chi'}+\sum_{i=1,..,s}
\epsilon^{H_i,\psi_i}_{\chi,\chi'} P_i. \label{bundle}
\end{equation} If $\chi_{|H_i}=\psi_i^{n_i(\chi)}$, then
\begin{equation}
 d_\chi+d_{\chi'}=d_{\chi \chi'}+\sum_{i=1,..,s}
\lfloor \frac{n_i(\chi)+n_i(\chi')}{e_i}
 \rfloor.
\label{degree} \end{equation} Let $H$ be the subgroup of $G$
generated by the inertia subgroups $H_i$ and let $Q=G/H$. One sees
that there is an exact sequence of groups $$1\lra Q^* \lra G^*\lra
H^* \lra 1.$$ The generators $\psi _i$ of $H_i^*$ define
isomorphisms $H_i^*\cong \Z _{e_i}$ where $e_i:=|H_i|$. Therefore,
we have an induced injective homomorphism $$\varphi:
H^*\hookrightarrow \prod_{i=1,..,s} \Z_{e_i}$$ such that the
induced maps $\varphi _i:H^*\lra \Z _{e_i}$ are surjective. By
abuse of notation, we will also denote by $\varphi$ the induced
homomorphism $\varphi: G^* \lra \prod_{i=1,..,s} \Z_{e_i}$. We
will write $$\varphi (\chi )=(n_1(\chi ),...,n_s(\chi ))\ \ \
\forall \chi \in G^*.$$ Let $\mu(\chi)$ be the order of $\chi$. By
\cite{Pa} Proposition 2.1, $$ d_\chi=\sum_{i=1,..,s}\frac{
n_i(\chi)}{e_i}.$$ We will now analyze all possible inertia groups
$H$.

{\bf Case 1:} $s=4$, and $e:=e_i=2$. Then $H^* \subset \Z_2^4$.\\
Note that $H^* \ne \Z_2^4$ since $(1,0,0,0) \notin H^*$. Thus
$H^*\cong (\Z_2)^s$ with $1\leq s \leq 3$.

By Example 1, all of these possibilities occour.

{\bf Case 2:} $s=3$ and $e_1 \ge 3$, $e_2=e_3=2$.\\ There must be
a character $\chi$ with $\varphi(\chi)=(1,n_2,n_3)$, and so
$$d_\chi=\frac{1}{e_1}+\frac{n_2}{2}+\frac{n_3}{2}$$ which is not
an integer. Therefore this case is impossible.

{\bf Case 3:} $s=2$ and $e_1,e_2 \ge 3$.\\ Assume that $e_1 >
e_2$. Since $G^* \to \Z_{e_1}$ is surjective, there is $\chi\in
H^*$ with $\varphi(\chi)=(1,n_2)$. Then
$$d_\chi=\frac{1}{e_1}+\frac{n_2}{e_2} <1$$ which is impossible.
So we may assume that $e=e_1=e_2\geq 3$ and $H^* \subset \Z_e^2$.
Let $\varphi(\chi)=(n_1,n_2)$. One has $d_\chi
=\frac{n_1+n_2}{e}$. Thus $n_2=e-n_1$ for any $\chi \ne 1$.
Therefore, $H^*=\{(i,e-i)| 0 \le i \le e-1\} \cong \Z_e$. By
Example 1, all of these possibilities occour.

From the above discussion, it follows that:
\begin{prop}\label{cover}
 Let $\phi: C\lra E$ be a $G$-cover with $E$ an elliptic curve
and $\dim H^0(\omega _C ^{\ot 3})^G=4$. Then either $\phi$ is
ramified over $4$-points and the inertia group $H$ is isomorphic
to $(\Z _2)^s$ with $s\in \{ 1,2,3\}$ or $\phi$ is ramified over
$2$-points and the inertia group $H$ is isomorphic to $\Z_m$ with
$m\geq 3$. \end{prop}

\bigskip
\bigskip

\begin{minipage}{13cm}
\parbox[t]{6.5cm}{
Jungkai Alfred Chen\\
Department of Mathematics\\
National Taiwan University\\
No.1, Sec. 4, Roosevelt Road\\
Taipei, 106, Taiwan\\
jkchen@math.ntu.edu.tw } \hfill
\parbox[t]{5.5cm}{Christopher D. Hacon\\
Department of Mathematics\\
University of Utah\\
155 South 1400 East, JWB 233\\
Salt Lake City, Utah 84112-0090, USA\\
hacon@math.ucr.edu}
\end{minipage}
\end{document}